\documentclass[11pt]{amsart}
\usepackage{amsfonts,amssymb,amsthm,amsmath,amsxtra,amscd}
\usepackage[all]{xy}
\usepackage[dvips]{graphics}

\theoremstyle{plain}
\newtheorem{thm}{Theorem}[section]
\newtheorem{lem}[thm]{Lemma}
\newtheorem{prop}[thm]{Proposition}
\newtheorem{cor}[thm]{Corollary}
\newtheorem{claim}[thm]{Claim}

\theoremstyle{definition}

\theoremstyle{remark}
\newtheorem{eg}[thm]{Example}

\newtheorem{rmk}[thm]{Remark}

\numberwithin{equation}{section}

\newcommand{\C}{{\mathbb{C}}}

\renewcommand{\P}{{\mathbb{P}}}

\renewcommand{\O}{\mathcal O}                             

\renewcommand{\H}{\mathcal H}
\renewcommand{\L}{\mathcal L}

\newcommand{\I}{\mathcal I}

\newcommand{\M}{\mathcal M}

\renewcommand{\SS}{\mathcal S}
\newcommand{\QQ}{\mathcal Q}

\newcommand{\f} {\phi}

\newcommand{\e} {\eta}

\newcommand{\n} {\nu}

\newcommand{\p} {\pi}
\renewcommand{\r} {\rho}
\newcommand{\s} {\sigma}

\newcommand{\D} {\Delta}
\newcommand{\G} {\Gamma}
\renewcommand{\S} {\Sigma}

\def\.{\cdot}
\def\^{\widehat}
\def\~{\widetilde}
\def\o{\circ}

\def\rat{\dashrightarrow}

\def\({\left(}
\def\){\right)}

\DeclareMathOperator{\codim} {codim}

\DeclareMathOperator{\Dv} {Dv}

\DeclareMathOperator{\Supp} {Supp}

\begin{document}

\title{Bad loci of free linear systems}

\subjclass[2000]{Primary 14C20; Secondary 14J25}
\keywords{Linear system, Bertini's theorems, Stein factorization}

\author[Tommaso de Fernex]{Tommaso de Fernex}
\address{Department of Mathematics, University of Michigan,
East Hall, 525 East University Avenue, Ann Arbor, MI 48109-1109, USA}
\email{{\tt defernex@umich.edu}}

\author[Antonio Lanteri]{Antonio Lanteri}
\address{Dipartimento di Matematica ``F. Enriques'',
Universit\`a degli Studi di Milano, via Saldini 50, 20133 Milano, Italy}
\email{{\tt lanteri@mat.unimi.it}}

\begin{abstract}
The bad locus of a free linear system $\L$ on a normal
complex projective variety $X$ is defined as the
set $B(\L) \subseteq X$ of points that are not contained
in any irreducible and reduced member of $\L$.
In this paper we provide a geometric description of such locus
in terms of the morphism defined by $\L$.
In particular, assume that $\dim X \ge 2$ and $\L$
is the complete linear system
associated to an ample and spanned line bundle.
It is known that in this case $B(\L)$ is empty unless $X$ is a surface.
Then we prove that, when the latter occurs,
$B(\L)$ is not empty if and only if $\L$ defines a morphism onto
a two dimensional cone,
in which case $B(\L)$ is the inverse image of the vertex of the cone.
\end{abstract}

\maketitle

\section{Introduction and statement of the results}

Consider a linear system $\L$ on a normal complex projective variety $X$.
Bertini's Second Theorem gives a condition for the existence
of irreducible and reduced members of $\L$. Of course, if there is
one such divisor in $\L$, then the general one will satisfy
these conditions. It is then natural to inquire about
the points of $X$ where only reducible or not reduced divisors of $\L$ pass
through. These, by definition, are the {\it bad points} for $\L$.
The set they form is called the {\it bad locus} of $\L$, and
is denoted by $B(\L)$:
$$
B(\L) = \{ p \in X \mid
\text{every $D \in \L - p$ is reducible or not reduced} \}.
$$
If $\L$ is positive dimensional, then this locus is precisely
the complement in $X$ of the set swept out by irreducible
and reduced divisors of $\L$.

In spite of their basic definition, it seems that
bad loci where explicitly discussed for the first time
only recently, by Besana, Di Rocco and the second author,
in~\cite{BDRL}. Ibidem, general properties
and interesting examples are studied when $X$ is a
smooth variety of dimension $\ge 2$ and $\L$ is a free linear system
associated to an ample line bundle, with particular
emphasis to the case when the linear system is complete.
Among other things, it is shown for instance that
under these assumptions bad loci are finite sets
and appear only in dimension 2. Moreover,~\cite{BDRL}
contains the complete classification
of pairs $(S,L)$, where $S$ is a smooth complex projective
surface and $L$ is an ample and globally generated line bundle
such that $B(|L|)\ne\emptyset$, under the
further assumption that $L^2 \le 11$ or
$L^2 = p^2$ for a prime $p$. As it is observed in~\cite{BDRL},
the cases occurring in this classification
suggest that non-empty bad loci commonly occur on surfaces
that are covers of cones.

The purpose of this paper is to continue the study of
bad loci in a more general setting, considering
base-point free linear systems on normal projective varieties.
Our aim is to provide a description of the bad locus of a free linear
system in terms of the geometry of the associated morphism.
The main reason for allowing singularities to the variety relies on the
method of investigation, as we use reduction
via Stein factorization.

We investigate bad points in parallel with
what we call {\it locally-bad points}.
The latter are those points on $X$ for which every member
of $\L$ passing through one of them, say $p$, is locally reducible
or not reduced at that point $p$ (that is,
the element of the local ring $\O_{X,p}$, induced
by a local equation of any such divisor, is reducible).
We denote by $B_{loc}(\L)$
the {\it locally-bad locus} of $\L$:
$$
B_{loc}(\L) = \{ p \in X \mid
\text{every $D \in \L - p$ is locally reducible or not reduced at $p$} \}.
$$
As it is obvious from the definitions, one always has
$B_{loc}(\L) \subseteq B(\L)$.

Our main results are the following geometric
descriptions of bad loci and locally-bad loci.
Certain notation appearing in the statements below
will be defined in the main body of the paper.
Specifically, we refer the reader to Section~2 for the definitions
of divisorial part $\Dv(V)$ and codimensional supports
$\Supp^1(V)$ and $\Supp^{\ge 2}(V)$ of a proper closed subscheme
$V \subset X$, and to Section~3 for the definition
of locally-bad locus $B_{loc}(D)$ of an effective Weil divisor $D$ on $X$.

\begin{thm}\label{main-bad}
Let $\L$ be a free linear system on a normal complex
projective variety $X$, and let $f : X \to \P^m$ be the associated
morphism.
Let $Y = f(X)$, and let $\H$ be the trace of $|\O_{\P^m}(1)|$ on $Y$.
Whenever $\dim Y \ge 2$, let
$$
E = \sum_{y \in Y} \Dv(f^{-1}(y)),
$$
and set $\S = f(E)$ ($\S$ is a finite set, and
$E$ is a Weil divisor on $X$; see Lemma~\ref{lem} below). Then
$B(\L) = \emptyset$ unless exactly one of the following occurs:
\begin{enumerate}
\item
$Y$ is a curve and $(X,\L) \ne (\P^1,|\O_{\P^1}(1)|)$, in which case
either
$(Y,\H) = (\P^1,|\O_{\P^1}(1)|)$
and $B(\L)$ is the locus swept out by those fibers of $f$ that
are reducible or not reduced, or
$(Y,\H) \ne (\P^1,|\O_{\P^1}(1)|)$ and $B(\L) = X$.
\item
$(Y,\H) = (\P^2,|\O_{\P^2}(1)|)$.
\item
$Y \subset \P^m$ is a two dimensional cone of degree $\ge 2$
and $B(\L) = \Supp\(f^{-1}(\S \cup \{ q \})\)$, where $q$
is the vertex of the cone.
\item
$(Y,\H)$ is not as in any of the previous cases,
and $B(\L) = \Supp\(f^{-1}(\S)\)$.
\end{enumerate}
Moreover, case~(b) does not occur if $\L$ is a complete
linear system.
\end{thm}

\begin{thm}\label{main-reduction}
With the assumptions and notation as in Theorem~\ref{main-bad},
let
$$
\D = \bigcup_{y \in Y} \Supp^1(f^{-1}(y))
\cap \Supp^{\ge 2}\(f^{-1}(y)\)
$$
whenever $\dim Y \ge 2$.
Then $B_{loc}(\L) = \emptyset$ unless
exactly one of the following occurs:
\begin{enumerate}
\item
$Y$ is a curve and $B_{loc}(\L)$ is the union of the locally-bad
loci of the fibers of $f$.
\item
$(Y,\H) = (\P^2,|\O_{\P^2}(1)|)$.
\item
$Y \subset \P^m$ is a two dimensional cone of degree $\ge 2$
and $B_{loc}(\L) = \Supp\(f^{-1}(q)\) \cup B_{loc}(E) \cup \D$,
where $q$ is the vertex of the cone.
\item
$(Y,\H)$ is not as in any of the previous cases,
and $B_{loc}(\L) = B_{loc}(E) \cup \D$.
\end{enumerate}
Moreover, case~(b) does not occur if $\L$ is a complete
linear system.
\end{thm}

When $\dim X \ge 2$ and $\L$ is an ample linear system
with non-empty bad locus, then we know from~\cite{BDRL} that
$X$ is a surface and $B(\L)$ is a finite set.
In particular, cases~(a) and~(d) of Theorem~\ref{main-bad}
do not occur. A fortiori, similar conclusions hold
for the locally-bad locus of $\L$.
Moreover, if we restrict to the case of complete linear systems
associated to ample and spanned line bundles,
we are then able to formalize the intuition
from~\cite{BDRL} for which bad points are
related to morphisms over cones, obtaining the following result:

\begin{cor}\label{ample-complete}
With the notation as in Theorem~\ref{main-bad}, assume that
$\dim X \ge 2$ and that $\L$ is the complete linear system associated
to an ample line bundle on $X$.
Then $B(\L) = B_{loc}(\L) = \emptyset$ unless $\dim X=2$ and
$Y \subset \P^m$ is a cone of degree $\ge 2$, in which case
$B(\L) = B_{loc}(\L) = \Supp(f^{-1}(q))$
where $q$ is the vertex of the cone.
\end{cor}

Now a few words on the organization of the paper.
After fixing some notation and discussing
the definition and basic properties of locally-bad
loci, we present in Section~4
a criterion for two dimensional cones. This,
although very simple, plays an important role
in the study of bad loci and locally-bad loci.
The subsequent section is devoted to the study
of bad loci and locally-bad loci of ample and
free linear systems. This
gives a result of some interest in itself,
stated in Theorem~\ref{ample}. The general set-up is then
considered in Section~6, where we split the problem into two parts
by considering the Stein factorization of the
morphism defined by the linear system under investigation.
Specifically, the proof of Theorems~\ref{main-bad} and~\ref{main-reduction}
is reduced to the case of ample and free linear systems
by analyzing how (locally-)bad loci change by pull-back
via morphisms with connected fibers.
The case of complete linear systems is addressed in the
last section of the paper, where we show that
the case when $(Y,\H) = (\P^2,|\O_{\P^2}(1)|)$,
although effective in general
(see Examples~\ref{eg0}--\ref{eg2} below),
does not occur if $\L$ is complete.

The first author has been partially supported by
the University of Michigan Rackham Research Grant and
Summer Fellowship.
Both authors have been partially supported by the MIUR of
the Italian Government in the framework of the National Research
Project ``Geometry on Algebraic Varieties" (Cofin 2002).
The authors would like to thank Gian Mario Besana for fruitful conversations and
the referee for useful comments and suggestions.

\section{Notation}

In this paper all varieties are defined
over the complex field and are implicitly
assumed to be positive dimensional,
irreducible, and reduced.

Let $X$ be a projective variety.
Unless otherwise specified,
by {\it linear system} on $X$
we will mean a linear subsystem, say $\L$,
of the complete linear system associated to
some line bundle $L$ on $X$.
We will also say that $\L$ is ample if $L$ is so.
For a point $p \in X$, we write $\L  - p$
to mean the linear subsystem of $\L$
of elements passing through $p$.
If $\L$ is very ample,
the {\it degree} of the pair $(X,\L)$
is the degree of $X$ in the
projective embedding given by $\L$.
We additionally say that $(X,\L)$ is a
{\it cone} if, in such embedding, $X$
is a cone of degree $\ge 2$.
With this definition,
we exclude the degenerate case of $(\P^n,|\O_{\P^n}(1)|)$.

Let $X$ be a normal projective variety.
Let $V \subset X$ be a proper closed subscheme, and
let $\I \subset\O_X$ denote its ideal sheaf.
Let
$$
\I = \bigcap \QQ_k
$$
be a primary decomposition of $\I$.
By ordering the $\QQ_k$ according to their
codimension in $X$, we can write
$$
\I = (\I)_{\codim=1} \cap (\I)_{\codim\ge 2}
$$
where $(\I)_{\codim=1}$ defines a subscheme of pure
codimension 1 and
$(\I)_{\codim\ge 2}$ defines one of codimension $\ge 2$.
Associated to the subscheme defined by
$(\I)_{\codim=1}$ we obtain a Weil divisor on $X$,
that we denote by $\Dv(V)$ and call {\it divisorial part}
of $V$. Explicitly,
$$
\Dv(V) = \sum m_i V_i,
$$
where the $V_i$ are the irreducible components of $V$
of codimension 1 in $X$ and, for each $i$,
$m_i$ is the multiplicity of $V$ at the
generic point of $V_i$.
We denote by
$$
\Supp^1(V) \quad\text{and}\quad \Supp^{\ge 2}(V)
$$
the supports of the schemes defined by
$(\I)_{\codim =1}$ and $(\I)_{\codim\ge 2}$.
Of course $\Supp^1(V) = \Supp(\Dv(V))$. We remark that,
although in general the primary components
of $(\I)_{\codim\ge 2}$ are not uniquely determined,
their associated primes are (see for instance~\cite{AM}, Theorem~4.5),
so $\Supp^{\ge 2}(V)$ is well defined.

Consider now a morphism $f : X \to \P^m$,
where $X$ is a normal, projective variety.
For any subscheme $V \subseteq X$,
we use $f(V)$ to denote the set-theoretic image of $V$,
and for a Weil divisor $D$ on $X$,
we put $f(D) := f(\Supp(D))$.
For a subscheme $W \subseteq f(X)$, $f^{-1}(W)$
will instead denote the scheme-theoretic inverse
image of $W$. In particular, for $q \in f(X)$, $f^{-1}(q)$
denotes the scheme-theoretic fiber over $q$.

\begin{lem}\label{lem}
If $\dim(f(X)) \ge 2$, then
the set $\S = \{q \in f(X) \mid \Dv(f^{-1}(q)) \ne 0 \}$ is finite.
\end{lem}

\begin{proof}
This is clear when $\dim X = 2$, as in this case
$f$ is a generically finite morphism. Suppose then that
$\dim X \ge 3$. By the semi-continuity of the dimension
of the fibers of $f$, $\S$ is a closed subset of $f(X)$.
Take a resolution of singularities $\~X \to X$, and let
$g : \~X \to \P^m$ be the composite morphism.
Fix a general hyperplane $H$ of $\P^m$, so that
$g^{-1}(H)$ is smooth. On the other hand $g^{-1}(H)$ is connected
since $g(\~X) = f(X)$ is at least two dimensional (\cite{Jou}, Theorem~7.1;
see also~\cite{FL}, Theorem~1.1).
Therefore $g^{-1}(H)$ is irreducible, and so is $f^{-1}(H)$.
This implies that $\S$ is zero dimensional.
\end{proof}

\section{Locally-bad points}

Let $X$ be a normal complex projective variety.
Given an effective Weil divisor $D$ on $X$ and a point $p \in D$,
we say that $p$ is a {\it locally-bad point} of
$D$ if, writing $D = \sum D_i$ with $D_i$ prime effective divisors on $X$
(repetitions are allowed),
$$
p \in \bigcup_{i \ne j} (D_i \cap D_j).
$$
We denote by $B_{loc}(D)$ the set of locally-bad points of $D$.
Rephrasing the definition, $B_{loc}(D)$ is the support of the non-reduced
components of $D$ union all intersections
of pairs of distinct irreducible components of $D$.

\begin{rmk}
The definition of locally-bad point has of course local
nature. When $D$ is a Cartier divisor on $X$, the locally-bad points of $D$
can be characterized as follow. Let $h \in \O_{X,p}$ be the
element induced by a local equation defining $D$ in a
neighborhood of a point $p \in D$.
Then $p \in B_{loc}(D)$ if and only if $h$ is reducible
in $\O_{X,p}$. Note also that, if $D$ is a Cartier divisor,
then $B_{loc}(D) = B_{loc}(\{D\})$, where
$\{D\}$ denotes the linear system on $X$ consisting
of the single element $D$. Moreover, we can rewrite
the definition of locally bad locus of a
linear system $\L$ on $X$, given in the introduction,
as follows:
$$
B_{loc}(\L) = \{ p \in X \mid
\text{$p \in B_{loc}(D)$ for every $D \in \L - p$} \} \\
$$
\end{rmk}

\begin{prop}\label{properties}
The locally-bad locus of a free linear system $\L$ is a proper
closed subset of $X$.
\end{prop}

\begin{proof}
Consider the incidence variety
$$
I = \{ (x,D) \in X \times \L \mid x \in D \}
$$
and the two projections $\p : I \to X$ and $\r : I \to \L$.
If $m = \dim \L$, then, for every $x \in X$, $\dim(\p^{-1}(x)) = m-1$
since $\L$ is free. Let
$\L_B = \{ D \in \L \mid B_{loc}(D) \ne \emptyset  \}$, and
consider the set
$$
I_B = \{ (x,D) \mid \text{$D \in \L_B$ and $x \in B_{loc}(D)$} \}.
$$
We claim that $I_B$ is closed in $I$.
To see this, consider a small pointed disk $0 \in T \subset \C$
parameterizing a non-trivial family of pointed divisors $(D_t,x_t)$,
such that $D_t \in \L$ for all $t \in T$ and $x_t \in B_{loc}(D_t)$
for $t \ne 0$. By the definition of locally bad point, and
possibly shrinking $T$ nearby $0$,  we can find families
of effective Weil divisors $A_t$ and $B_t$ on $X$,
parameterized by $T \setminus \{0\}$, such that $D_t = A_t + B_t$
and $x_t \in A_t \cap B_t$. Then $A_t$ and $B_t$ degenerate, respectively,
to two Weil divisors $A_0$ and $B_0$, and $D_0 = A_0 + B_0$.
We are not excluding that either of $A_t$ and $B_t$
might have components that are fixed throughout the deformation,
nor that they might share common components at any stage of the
degeneration.
Since $x_0 \in A_0 \cap B_0$, we conclude that $x_0 \in B_{loc}(D_0)$.
This, combined with the fact
that $B_{loc}(D)$ is closed in $D$ for every $D \in \L$,
implies that $I_B$ is closed in $I$, as claimed.
Observing that $B_{loc}(\L)$ is the locus of the fibers
of (maximal) dimension $m-1$ of $\p|_{I_B}$,
we conclude that it is
closed in $X$ by the semi-continuity of the fiber-dimension.
Since $\L$ is free and $X$ is normal, the generic element
of $\L$ is smooth in codimension 1 by Bertini Theorem;
in particular, it has empty locally-bad locus. This implies that
$\dim \L_B < \dim \L = m$.
We thus see that $\dim I_B < \dim I = m + \dim X -1$.
But $I_B = \p^{-1}(\p(I_B))$ and $B_{loc}(\L) \subset \p(I_B)$, hence
$\dim B_{loc}(\L) < \dim X$.
\end{proof}

\section{Bad loci of very ample linear systems}

As it is probably well known,
two dimensional cones can be characterized in terms of
the reducibility of their hyperplane sections:

\begin{prop}\label{cones}
Let $Y \subset \P^m$ be a projective variety
of dimension $\ge 2$, and let
$q \in Y$. Then every hyperplane section $H \subset Y$ passing through $q$
is reducible or not reduced at $q$
if and only if $\dim Y = 2$ and
$Y$ is a cone of degree $\ge 2$, with vertex $q$.
\end{prop}

\begin{proof}
Let $\H$ be the trace of $|\O_{\P^m}(1)|$ on $Y$.
The ``if'' part of the statement is clear: if $(Y,\H)$ is a
two dimensional cone with vertex $q$, then every
hyperplane section $H$ through $q$ is the cone over a number of points
equal to the degree of $Y$ (counted with multiplicities), thus
$H$ fails to be irreducible and reduced.
Conversely, assume that $(Y,\H)$ is not a
two dimensional cone with vertex $q$. Let
$\p : \P^m \rat \P^{m-1}$ be the linear projection
from $q$, and set $\f = \p|_Y$. Then $\f$ is precisely the
rational map defined by $\H - q$. Since
$Y$ is not a two dimensional cone with vertex $q$, we get $\dim \f(Y) \ge 2$.
Then Bertini's Second Theorem (see for instance~\cite{Kle}) implies that
the general element of $\H - q$ is irreducible and reduced.
\end{proof}

From Proposition~\ref{cones} we immediately obtain
the following criterion for bad points of very ample linear systems.

\begin{cor}\label{very-ample}
Let $Y$ be a projective variety
and let $\H$ be a very ample linear system on $Y$. Let
$q \in Y$. Then the following are equivalent:
\begin{enumerate}
\item   $q$ is a bad point for $\H$.
\item   $q$ is a locally-bad point for $\H$.
\item   $(Y,\H)$ is a two dimensional cone with vertex $q$.
\end{enumerate}
\end{cor}

\section{Bad loci of ample and free linear systems}

In this section we study bad loci and locally-bad loci of
ample and free linear systems.

\begin{thm}\label{ample}
Let $\L$ be an ample and free linear system on a normal complex projective
variety $X$. Denote by $f : X \to \P^m$
the morphism defined by $\L$, let $Y = f(X)$, and let
$\H$ be the trace of $|\O_{\P^m}(1)|$ on $Y$. Then
$B(\L) = B_{loc}(\L) = \emptyset$ unless
one of the following cases occurs:
\begin{enumerate}
\item
$X$ is a curve and $(X,\L) \ne (\P^1,|\O_{\P^1}(1)|)$, $B(\L) = X$,
and $B_{loc}(\L)$ is the ramification locus of $f$.
\item
$(Y,\H) = (\P^2,|\O_{\P^2}(1)|)$.
\item
$(Y,\H)$ is a two dimensional cone
and $B_{loc}(\L) = B(\L) = \Supp\(f^{-1}(q)\)$,
where $q$ is the vertex of the cone.
\end{enumerate}
\end{thm}

\begin{proof}
Since $\L$ is ample, we know by~\cite{BDRL}, Theorem~2,
that $B(\L) = \emptyset$ whenever $\dim X > 2$.
(Although $X$ is therein assumed to be smooth,
the proof of Theorem~2(i) in~\cite{BDRL} is still
valid for any complex projective variety.) Since
$B_{loc}(\L) \subseteq B(\L)$, the same conclusion holds for $B_{loc}(\L)$.
Therefore we can assume that $\dim X \le 2$.

Assume that $X$ is a (smooth) curve. Clearly $B(\L) = X$ unless
$(X,\L)$ has degree 1, in which case
$(X,\L) = (\P^1,|\O_{\P^1}(1)|)$ and
$B(\L) = B_{loc}(\L) = \emptyset$.
It is also clear that the locally-bad points of the fibers
of $f$ are precisely those points where the
differential of $f$ vanishes. This gives the
description of $B_{loc}(\L)$ when $X$ is a curve.

For the remainder of the proof, we will assume that
$\dim X = 2$.
We also suppose that $(Y,\H) \ne (\P^2,|\O_{\P^2}(1)|)$.
Fix a point $p$ in $X$ and set $q = f(p)$.

We need to prove that $p$ is a bad point for $\L$ if any only if
$p$ is a locally-bad point for $\L$, and this if and only if
$Y$ is a cone with vertex $q$. The first condition
is of course implied by the second one. The second condition, in turns, is
implied by the third one: if $Y$ is a cone with vertex $q$,
then $B_{loc}(\L)$ contains the
support of the fiber over $q$. Indeed,
by Proposition~\ref{cones}, $q$ is a locally-bad point of
every hyperplane section $H$ of $Y$
passing through $q$, and so $p$ is a locally-bad point of the
pull back $f^*H$. Therefore, to close the circle of
implications and conclude the proof of the theorem, we are left with
showing the following

\begin{claim}\label{claim}
Let $(Y,\H) \ne (\P^2,|\O_{\P^2}(1)|)$, let $q \in Y$, and
suppose that $(Y,\H)$ is not a cone with vertex $q$.
Then, for a general $H \in \H - q$, $f^*H$
is irreducible and reduced.
\end{claim}

Since $(Y,\H)$ is not a cone with vertex $q$,
we can find a regular point $q' \ne q$ of $Y$ such that
$f$ is locally \'etale over $q'$ and the line
in $\P^m$ joining $q$ and $q'$ is not contained in $Y$.
Let $p_1,\dots, p_d$ be the pre-images of $q'$.
Consider the pencil $\H'$ generated by two general elements
$H_1, H_2 \in \H - q - q'$, and let $\L' = f^*\H'$.
Note that $\L'$ is a pencil without fixed components,
and defines a rational map $\f : X\rat \P^1$.
Let
$$
\xymatrix{
\~ X \ar[d]_{\s} \ar[dr]^{\~ \f} \\
X \ar@{-->}[r]_{\f} & \P^1
}
$$
be a resolution of indeterminacy of $\f$ such that $\~X$
is a smooth surface.
Clearly, every fiber of $\~ \f$ is reducible or not reduced
if and only if $f^*H$ is reducible or not reduced for every $H \in \H'$.
Since $\~ X$ is a smooth surface, the general fiber
of $\~ \f$ is smooth. Therefore, if we suppose that the statement of the
claim is false, then
the general fiber of $\~ \f$ must be not connected. Thus, it suffices
to show that such situation leads to a contradiction.

Take the Stein factorization of $\~ \f$:
$$
\xymatrix{
\~ X \ar[dr]_{\~ \f}  \ar[r]^{\e} & W \ar[d]^{\n} \\
& \P^1.
}
$$
$W$ is a normal (hence smooth) curve,
$\e$ is a morphism with connected fibers
and $\n$ is a finite morphism of degree $\geq 2$.
It is clear that the points of $\P^1$ parameterize the elements of
$\L'$ and the points of $W$ parameterize components (not necessarily
irreducible) of such divisors.
For $i = 1,\dots,d$, we define $V_i \subset W$ to be the set of points $w$
which correspond to components $C_w$ (of divisors in $\L'$)
passing through $p_i$:
$$
V_i = \{ w \in W \,|\, C_w \ni p_i \}.
$$
These sets are closed in $W$.
Let $U \subset \P^1$ be the set parameterizing curves $D \in \L'$
whose image $H = f(D)$ is irreducible.
We have, from Proposition~\ref{cones}, that the general hyperplane
section of $Y$ is irreducible. Therefore, by generality of the choice
of the pencil $\H'$, we conclude that $U$ is a non-empty open set of $\P^1$,
so that $\n^{-1}(U)$ is an open set of $W$.
Since
$\{p_1,\dots,p_d \} = f^{-1}(q')$ and $q'$ is a base point of $\H'$,
the sets $V_i$, for $i=1,\dots,d$, cover $\n^{-1}(U)$:
$$
\bigcup_{i=1}^d V_i \supseteq \n^{-1}(U).
$$
But $\n^{-1}(U)$ is open and the $V_i$ are closed, thus
for at least one index $i$, say $i=1$, $V_1 \supseteq \n^{-1}(U)$.
Hence $V_1 = W$.
Note that, if $u \in \P^1$ is a sufficiently general point
and $D_u$ is the
divisor corresponding to $u$, then $\n^{-1}(u)$ contains at least
two distinct points $w,w'$ of $W$. Since $W = V_1$, we conclude
that the two components of $D_u$ corresponding to $w$ and $w'$
pass through $p_1$. But this is impossible.
Indeed, since $Y$ is smooth at $q'$ and
$f$ is locally \'etale over $q'$,
two distinct components of a general divisor $D \in \L'$ cannot meet
at a point lying over $q'$.
This proves the claim, hence concludes the proof of the theorem.
\end{proof}

\section{General setup and reduction via Stein factorization}

Let $X$ be a normal complex projective variety and
$\L$ be a free linear system on $X$. Let
$m = \dim \L$, and denote by
$f : X \to \P^m$
the morphism defined by $\L$.
Let $Y = f(X)$ and $\H$ be the trace of $|\O_{\P^m}(1)|$ on $Y$,
so that $\L = f^*\H$.
Let
$$
\xymatrix{
X \ar[r]^g \ar[rd]_f & Z \ar[d]^h \\
& Y
}
$$
be the Stein factorization of $f$, and
set $\M = h^*\H$.
Recall that $g$ has connected fibers, $h$ is a finite morphism,
and $Z$ is a normal variety. In particular, $\M$ is ample.

Let $y$ be a point in $Y$. For every $z \in Z$ with $h(z) = y$,
we denote by $F_z$ the connected component of $f^{-1}(y)$
that is mapped to $z$, so that
$$
f^{-1}(y) = \bigsqcup_{h(z)=y} F_z.
$$

\begin{rmk}\label{F_z}
Note that $F_z = g^{-1}(z)$ whenever $h$ is
a local isomorphism near $z$.
\end{rmk}

If $\dim Z=1$, we denote by $F$ the sum of all those $F_z$,
as $z$ varies in $Z$, that are reducible or non-reduced:
$$
F = \sum_{\substack{\text{$F_z$ is reducible}\\ \text{or not reduced}}}F_z.
$$
Note that $F$ is a Weil divisor; indeed, $F_z$ is smooth
for a general $z \in Z$, since $X$ is normal, and so
it is irreducible and reduced because connected.
If $\dim Z \ge 2$, let
$$
\D = \bigcup_{q \in Y} \Supp^1(f^{-1}(q)) \cap
\Supp^{\ge 2}\(f^{-1}(q)\)
$$
and
$$
E = \sum_{q \in Y} \Dv(f^{-1}(q)).
$$
Note that $E$ is a Weil divisor, by Lemma~\ref{lem}.

With the following three propositions, we relate the bad locus
of $\L$ to that of $\M$.
This will allow us to reduce the
proof of Theorems~\ref{main-bad} and~\ref{main-reduction} to the
case when the linear system is ample, that was settled in Theorem~\ref{ample}.

\begin{prop}\label{bad}
Consider the notation introduced at the beginning of the section.
If $\dim Z = 1$,
then $B(\L) = X$ unless $(Z,\M) = (\P^1,|\O_{\P^1}(1)|)$,
in which case $B(\L) = F$.
If $\dim Z > 1$, then $B(\L) = \Supp\(g^{-1}(B(\M) \cup g(E))\)$.
\end{prop}

\begin{proof}
If $\dim Z = 1$, then $Z$ is a smooth curve,
and every element of $\M$ is reducible or non-reduced
unless $(Z,\M) = (\P^1,|\O_{\P^1}(1)|)$.
In this event, since the elements of $\L$ are precisely
the fibers of $g$, the first part follows at once.
Assume now that $\dim Z \ge 2$. Clearly
$g^{-1}(B(\M)) \subseteq B(\L)$. Now we consider
a point $z \in Z \setminus B(\M)$.
Note that $\M - z$ has finite base locus, since $\M$
is ample and free.
Thus, for any fixed $p \in g^{-1}(z)$, we have $p \in B(\L)$
if and only if $z \in g(E)$.
\end{proof}

\begin{prop}\label{reduction-dim=1}
With the notation as introduced at the beginning of the section,
assume that $\dim Z = 1$.
Then $B_{loc}(\L) = \Supp\(g^{-1}(B_{loc}(\M))\) \cup B_{loc}(F)$.
\end{prop}

\begin{proof}
If $\dim Z = 1$, then $Z$ is smooth and every fiber of $g$ is a Cartier
divisor on $X$. Clearly the support of $g^{-1}(B_{loc}(\M))$
is included in $B_{loc}(\L)$.
Indeed, if $z \in B_{loc}(\M)$ and $M \in \M-z$, then $g^*M$ is
non-reduced along $g^{-1}(z)$ because $M$
is singular at $z$, hence
each point of $g^{-1}(z)$ is a locally bad point of $g^*M$.

Suppose then that $p \not \in g^{-1}(B_{loc}(\M))$, and let $z = g(p)$.
Since the locally-bad points of $\M$
are precisely the ramification points of $h$, we have
$g^{-1}(z) = F_z$ (see Remark~\ref{F_z}).
By our choice of $p$, the general member of $\L - p$
consists of $F_z$ plus, possibly, some other divisors disjoint from $F_z$.
Therefore $p \in B_{loc}(\L)$ if and only if $p \in B_{loc}(F_z)$.
Since $B_{loc}(F) = \cup_{z \in Z}B_{loc}(F_z)$, the conclusion follows.
\end{proof}

\begin{prop}\label{reduction-dim>1}
With the notation introduced at the beginning of the section,
assume that $\dim Z \ge 2$ and $(Y,\H) \ne (\P^2,|\O_{\P^2}(1)|)$.
Then
$B_{loc}(\L) = \Supp\(g^{-1}(B_{loc}(\M))\) \cup B_{loc}(E) \cup \D$.
\end{prop}

\begin{proof}
Let $z \in B_{loc}(\M)$ and $M \in \M-z$. We start observing that $g^*M$ is
non-reduced along the support of $\Dv(g^{-1}(z))$, since $M$
is singular at $z$.

\begin{lem}\label{M}
Counting multiplicities,
there are at least two components of $g^*M$ passing
through any given point of
$\Supp^{\ge 2}(g^{-1}(z)) \setminus \Supp^1(g^{-1}(z))$.
\end{lem}

\begin{proof}[Proof of Lemma~\ref{M}]
We use the description
of $(Y,\H)$, as it follows from Theorem~\ref{ample} applied to $(Z,\M)$.
Since we are assuming that $z$ is a locally-bad point of $\M$,
and in view of the assumptions of the proposition, $(Y,\H)$
is a two dimensional
cone with vertex $q = h(z)$.
We stratify $Y$ according to the codimension of the fibers of $f$:
$$
Y = \bigsqcup_{k=1}^n Y_k, \quad\text{where}\quad
Y_k = \{ y \in Y \mid \codim_X(f^{-1}(y)) = k \}.
$$
Here $n = \dim X - \dim Y$, $Y_1 = f(E)$, and $Y_n$
is a non-empty open subset of $Y$.

Let $\SS$ be the set of lines $\ell \subset Y$, passing through $q$,
that are not contained in the singular locus of $Y$
and intersect every $Y_k \setminus \{q\}$ properly.
Then $\SS$ is parameterized by a non-empty open subset of
the curve which $Y$ is the cone over. In particular, $\SS$ is connected.
For a line $\ell \in \SS$, we define $\ell' \subset X$ to be the closure
of the support of $f^{-1}(\ell \cap Y_n)$.
Note that, for every irreducible component $W$
of $\Supp^{\ge 2}(g^{-1}(z))$ that is not
contained in $\Supp^1(g^{-1}(z))$, the set
$$
\SS_W = \{ \ell \in \SS \mid \ell' \supset W \}
$$
is closed in $\SS$.

Since the property stated in the
lemma is closed in $\M - z$, we can assume that $M$ is
a general element in $\M - z$. Then $M = h^*H$, where
$H = \ell_1 + \dots + \ell_d$
is the sum of $d$ lines $\ell_i \subset Y$ each passing through $q$,
$d$ being the degree of $(Y,\H)$. Furthermore, we can assume that
each $\ell_i$ is a general element of $\SS$.
In particular, we have that
$$
f^*H = \ell_1' + \dots + \ell_d' + \Dv(f^{-1}(q)),
$$
where $\ell_i'$ is the closure
of the support of $f^{-1}(\ell_i \cap Y_n)$.
Since $W$ is contained in $f^*H$ but is not contained in $\Dv(f^{-1}(q))$,
some of the $\ell_i'$'s must contain $W$. Since we are choosing
$H$ general (recall also that $\SS_W$ is closed in $\SS$),
this is possible if and only if
$\SS_W = \SS$. But this precisely means that
every $\ell_i'$ contains $W$. Since $d \ge 2$, this proves the lemma.
\end{proof}

We now continue the proof of Proposition~\ref{reduction-dim>1}.
Lemma~\ref{M}, combined with the starting observation, implies that
each point of $g^{-1}(z)$ is a locally bad point of $g^*M$.
In particular, this implies that
$$
\Supp\(g^{-1}(B_{loc}(\M))\) \subseteq B_{loc}(\L).
$$
It is clear that $B_{loc}(E) \subseteq B_{loc}(\L)$.
Therefore we fix, for the remainder of the proof, a point
$$
p \in X \setminus \left(\Supp\(g^{-1}(B_{loc}(\M))\)
\cup B_{loc}(E)\right).
$$
To complete the proof, we need to show that $p$ is a
locally-bad point for $\L$ if and only if $p \in \D$.

We set $z = g(p)$ and $q = h(z)$.
Since the linear forms on $\P^m$ vanishing at $q$
generate the maximal ideal of $q$ in $Y$,
the base scheme of $\H - q$ is precisely the
reduced point $q$. Therefore
the base scheme of $f^*(\H-q)$ is
the fiber $f^{-1}(q)$. Thus for any $H \in \H-q$
we can write
$$
f^*H = \Dv(f^{-1}(q)) + A.
$$
Moving $H$ in $\H-q$, $A$ moves in a linear system
of Weil divisors whose base locus is precisely
$\Supp^{\ge 2}(f^{-1}(q))$.

Since $X$ is normal and $g$ is a morphism with connected fibers,
the generic fiber of $g$ is
irreducible and reduced. Thus we can find a non-empty open subset
$U \subset Z$ over which $g$ is equidimensional (of dimension
$\dim X - \dim Z$) and has irreducible and reduced fibers.
Then, if $H \in \H - q$ is sufficiently general and $M = h^*H$,
we can assume that the following conditions are satisfied:
\begin{enumerate}
\item
$z \not \in B_{loc}(M)$ (recall that $z \not \in B_{loc}(\M)$), and
\item
every irreducible component of $M$ intersects $U$.
\end{enumerate}
Note indeed that, to satisfy~(b), it is enough to assume that each
irreducible component of $H$ intersects $Y \setminus h(Y \setminus U)$.
Condition~(b) implies that the fiber
of $g$ over the generic point of any given irreducible
component of $M$ is irreducible and reduced
and has dimension $\dim X - \dim Z$. Since each irreducible
component of $A$ maps onto an irreducible
component of $M$, condition~(a)
implies that $p \not \in B_{loc}(A)$.
Since $p \not \in B_{loc}(E)$,
we deduce that $p$ is a locally-bad point for $\L$
if and only if $p \in \Dv(f^{-1}(q)) \cap A$.
By varying $H$ among general members of $\H - q$, we conclude
that $p \in B_{loc}(\L)$ if and only if $p \in \D$,
as we needed.
\end{proof}

\section{Further remarks and the case of complete linear systems}

It appears from Theorems~\ref{main-bad}
and~\ref{main-reduction} that bad points
and locally-bad points
commonly arise from the divisorial components of the fibers
of the morphism defined by the linear system,
and as pre-images of vertexes of two dimensional
cones. Explicit examples of bad points
where $\L$ is a complete ample linear system on a smooth surface
and $Y$ is a cone are discussed in
Section~3 of~\cite{BDRL}.

There is a further case listed in Theorems~\ref{main-bad}
and~\ref{main-reduction}
that appears somehow anomalous, namely when $\L$ is
a linear system on a surface having non-trivial (locally-)bad locus
and $(Y,\H) = (\P^2,|\O_{\P^2}(1)|)$.
The following examples show that this case is effective.
The first one is obvious, while the next two
are more interesting, since they involve nets on smooth
surfaces.

\begin{eg}\label{eg0}
Let $X \subset \P^m$ be a two dimensional cone of degree $\ge 2$,
with vertex a point $p$. Let $\P^m \rat \P^2$ be any linear
projection from a center disjoint from $X$, and let
$f : X \to \P^2$ be the induced morphism. Then,
if $\L = f^*|\O_{\P^2}(1)|$, we have
$B(\L) = B_{loc}(\L) = \{p\}$.
\end{eg}

\begin{eg}\label{eg1}
Consider the following three situations:
\begin{enumerate}
\item[(i)]
Let $X = \P^2$, fix two lines $A_1$ and $A_2$, and let $x$ be the
point of intersection of these two lines. For any $d \ge 2$, let
$C \subset X$ be a general curve of degree $d$ and consider the
linear system $\L$ spanned by $C$, $dA_1$ and $dA_2$.
\item[(ii)]
Let $X$ be a Del Pezzo surface of degree $K_X^2=1$, let $x$ be the
base point of the anticanonical pencil $|-K_X|$, and fix two
elements $A_1,A_2 \in |-K_X|$ . For $d \geq 3$ the linear
system $|-dK_X|$ is very ample. Let $C \in |-dK_X|$ be a general
curve, and consider the linear system $\L$ spanned by $C$, $dA_1$
and $dA_2$.
\item[(iii)]
Let $X$ be a minimal surface of general type with $K_X^2=1$ and
$p_g(X)=2$, let $x$ be the base point of the canonical pencil
$|K_X|$ and fix two elements $A_1, A_2 \in |K_X|$. It is
well known that the linear system $|dK_X|$ is very ample for
$d \geq 5$ (see~\cite{Kod}, or~\cite{Bom}). For any such $d$, let
$C \in |dK_X|$ be a general curve and consider the linear system $\L$
spanned by $C$, $dA_1$ and $dA_2$.
\end{enumerate}
In all these cases $\L$ is a free linear system defining a
morphism of degree $d^2$ onto $\P^2$, and
$$
B(\L)=B_{loc}(\L)=\{x\}.
$$
This can be seen as follows. First note
that $B(\L)=B_{loc}(\L)$ as it follows from the proof of Theorem 2
in~\cite{BDRL}, and clearly this set contains the point $x$.
It is easy to see that every other
element of the pencil generated by $dA_1$ and $dA_2$ is the sum of
$d$ distinct irreducible curves, any two of which meet only at $x$, and
this pencil sweeps out the whole $X$. Hence, for any given point of
$X \setminus (A_1 \cup A_2)$, we can find an
element of the pencil having exactly one irreducible (and reduced)
component passing through it. In particular, the bad locus of this
pencil is contained in $A_1 \cup A_2$. Now, letting
$L$ be the line bundle associated with $\L$, we show how the
general choice of $C \in |L|$ can guarantee that the corresponding
linear subsystem $\L \subset |L|$, as defined in the examples, has
only $x$ as a bad point. Consider the locus $W \subseteq |L|$
defined by those curves $C$ such that the corresponding $\L$ has
some bad point different from $x$. Since $B(\L)$ is contained
in the bad locus of the pencil generated by $dA_1$ and $dA_2$,
$W$ is equal to the union of the loci $W_y$ consisting of those
$C \in |L|$ for which $y$ is a bad point of $\L$, for $y$ ranging
on $A_1 \cup A_2 \setminus \{x\}$. Note that $\dim W \leq \dim W_y +1$
for any such $y$. On the other hand, due to the very ampleness
of $L$, the chain of inclusions
$$
W_y \cap |L-y| \subseteq |L - 2y| \subset |L|
$$
implies that the codimension of $W_y$ in $|L|$ is at least $2$.
Therefore $W$ has positive codimension inside $|L|$, and
we can choose $C$ outside it.
\end{eg}

\begin{eg}\label{eg2}
For $d \ge 2$, consider the following situations:
\begin{enumerate}
\item[(j)]
Let $x \in \P^2$ and $\G_1, \G_2 \in
|\O_{\P^2}(d)-dx|$. Each $\G_i$ consists of $d$ lines (not
necessarily distinct) through $x$. Suppose that
$\Supp(\G_1) \cap \Supp(\G_2) =\{x\}$. Then the
linear system $\L$ spanned by $\G_1$, $\G_2$ and a general curve
of degree $d$ is free, and
$B(\L)= B_{loc}(\L)=\{x\}$. This is slightly more general than
Example~\ref{eg1}(i).
\item[(jj)]
Now consider two distinct points $x,y \in \P^2$, and
let $\G_1 \in |\O_{\P^2}(d)-dx|$ and
$\G_2 \in |\O_{\P^2}(d)-dy|$ be general elements. We can assume that
neither $\G_1$ nor $\G_2$ contains the line $\ell$ joining $x$ and
$y$. Then the linear system $\L$ spanned by
$\G_1$, $\G_2$ and $d\ell$ is free, and $B(\L)= B_{loc}(\L)= \{x,
y\}$.
\item[(jjj)]
Finally, let $\ell_1, \ell_2, \ell_3$ be three lines of $\P^2$ not
in a pencil, and let $\L$ be the net spanned by $d\ell_1, d\ell_2,
d\ell_3$. Then $\L$ is a free linear system with three bad points,
which are also locally-bad points.
\end{enumerate}
The assertion in~(j) follows from the same arguments used in
Example~\ref{eg1}, once we observe that the bad locus of the pencil
spanned by $\G_1$ and $\G_2$ is contained in a curve. The
assertions in~(jj) and~(jjj) can be directly checked by
computations in appropriate affine coordinates.
\end{eg}

The case in which $(Y,\H) = (\P^2,|\O_{\P^2}(1)|)$
might be not too surprising, since the latter
may actually be viewed as a degenerate ``cone'' of degree 1.
However Example~\ref{eg2} suggests that this case
should be considered as a special one.
Indeed the existence of two (three) bad points mapped to two (three)
distinct points on $\P^2$ shows that it is impossible to
choose, in a meaningful way, a ``vertex'' for $(\P^2,|\O_{\P^2}(1)|)$.
It turns out in fact that this peculiar case
does not occur when $\L$ is a complete linear system.

\begin{prop}\label{complete}
Under the assumptions of Theorem~\ref{main-bad}, assume that
$\L$ is a complete linear system.
If there is a bad point for $\L$, then
$(Y,\H) \ne (\P^2,|\O_{\P^2}(1)|)$.
\end{prop}

\begin{proof}
If $\dim Y \ne 2$, there is nothing to prove, so we may assume that
$\dim Y = 2$. In particular, $\dim \L \ge 2$.
Then it is enough to show that, if there is a bad point for $\L$,
then we actually have $\dim \L \ge 3$. Suppose that $p$
is a bad point for $\L$. An arbitrary element $D \in \L - p$
can be written in the form $D = A + B$, where $A$ and $B$
are two effective divisors. Note that
both $A$ and $B$ span positive dimensional
linear systems on $X$, since the base locus of $\L-p$ is
finite by the ampleness of $\L$. Therefore
$\dim |A| \ge 1$ and $\dim |B| \ge 1$.
From the multiplication map
$$
H^0(\O_X(A)) \times H^0(\O_X(B)) \to H^0(\O_X(A+B)),
$$
we have a map
$$
|A| \times |B| \to |A+B| = \L.
$$
This map is clearly finite. On the other hand, there are elements of
$\L$ that are irreducible and reduced, by Bertini's Second Theorem.
Therefore the above map is not surjective, hence $\dim \L \ge 3$.
\end{proof}

In closing, a few indications on how to combine the results established
in Sections~5--7
to prove the theorems stated in the introduction.
We come back to the set-up of Section~5. For the proof of
Theorem~\ref{main-bad}, we apply Theorem~\ref{ample} to
$(Z,\M)$ and then use Propositions~\ref{bad};
the last claim of Theorems~\ref{main-bad}
follows from Proposition~\ref{complete}.
Similarly, we apply Theorem~\ref{ample},
Propositions~\ref{reduction-dim=1} and~\ref{reduction-dim>1},
and Proposition~\ref{complete}
for the proof of Theorem~\ref{main-reduction}.
Finally, combining Theorem~\ref{ample} and Proposition~\ref{complete},
we get Corollary~\ref{ample-complete}.

\providecommand{\bysame}{\leavevmode \hbox \o3em
{\hrulefill}\thinspace}

\end{document}